\documentclass[11pt,a4paper]{article}
\usepackage{latexsym}
\usepackage{graphicx}
\usepackage{amsmath}
\usepackage{amssymb}
\usepackage{epsfig}

\makeatletter \@addtoreset{equation}{section}
 \makeatother

\begin{document}

\title{ \textbf{$\gamma^*$-Regular,$\gamma$-Locally Compact and $\gamma$-Normal spaces}}
\author{$Bashir Ahmad^*$ and $Sabir Hussain^{**}$\\
{*} Centre for Advanced Studies in Pure and Applied Mathematics,
 Bahauddin Zakariya \\ University, Multan, Pakistan.
{\bf Present Address:} Department of Mathematics, King \\ Abdul
Aziz University P.O.Box 80203, Jeddah 21589, Saudi Arabia.\\
E.mail: drbashir9@gmail.com
\\{**} Department of Mathematics, Islamia University Bahawalpur,
Pakistan.\\ {\bf Present Address:} Department of Mathematics,
Yanbu University,\\ P. O. Box 31387, Yanbu Alsinaiyah, Saudi
Arabia.\\E. mail: sabiriub@yahoo.com.}
\date{}
\maketitle
\textbf{Abstract.} We continue studying the properties of
$\gamma_0$-compact, $\gamma^*$-regular and $\gamma$-normal spaces
defined in [5]. We also define and discuss $\gamma$-locally
compact spaces.
\\
AMS(2000)Subject Classification. Primary 54A05, 54A10, 54D10.
\section {Introduction}
S. Kasahara [8] introduced the concept of $\alpha$-closed graphs
of a function by  using the concept of  an operation $\alpha$ on
topological spaces in 1979. $\alpha$-closed sets were defined by
D.S. Jankovic [7] and he studied the functions with
$\alpha$-closed graphs. H. Ogata [9] introduced the notions of
$\gamma$-$T_i$ ,i = 0,1/2,1,2; $(\gamma,\beta)$- homeomorphism and
studied some topological properties.\\ B.Ahmad and F.U.Rehamn (
[10] respt. [2] ) defined and investigated several properties of
 $\gamma$-interior, $\gamma$-exterior, $\gamma$-closure and
$\gamma$-boundary points in topological  spaces. They also
discussed their properties in product spaces and studied the
characterizations of $(\gamma,\beta)$-continuous mappings
initiated by H. Ogata [9]. In 2003( respt. in 2005), B.Ahmad and
S.Hussain [3] (respt. [4]) continued studying the properties of
$\gamma$-operations on topological spaces introduced by S.
Kasahara [8]. They also defined and discuss several properties of
$\gamma$-nbd, $\gamma$-nbd base at x, $\gamma$-closed nbd,
$\gamma$-limit point, $\gamma$-isolated point, $\gamma$-convergent
point and $\gamma^*$-regular space. They further defined
$\gamma$-normal spaces, $\gamma_0$-compact, $\gamma^s$-regular and
 $\gamma^s$-normal in topological spaces ( [5] respt. [6] ) and
established many interesting properties.

        In this paper we continue to discuss the properties of
$\gamma_0$-compactness, $\gamma^*$-regular, $\gamma$-normal spaces
defined in [5]. We also define and discuss $\gamma$-locally
compact spaces.

        First, we recall preliminaries used in the sequel. Hereafter we shall write space
in place of topological space.\\
{\bf Definition [9].}  Let (X,$\tau$) be a space. An operation
$\gamma$ : $\tau\rightarrow$ P(X) is a function from $\tau$ to the
power set of X such that V $\subseteq V^\gamma$ , for each V
$\in\tau$, where $V^\gamma$ denotes the value of $\gamma$ at V.
The operations defined by $\gamma$(G) = G, $\gamma$(G) = cl(G) and
$\gamma$(G) = intcl(G) are examples of operation $\gamma$.\\
{\bf Definition [10].}  Let A $\subseteq$ X. A point a $\in$ A is
said to be $\gamma$-interior point of A iff there exists an open
nbd N of a such that $ N^\gamma\subseteq$ A and we denote the set
of all such points by $int_\gamma$(A). Thus
\begin {center}
        $int_\gamma (A) = \{ x \in A : x \in  N \in \tau$ and $N^\gamma\subseteq A \}$.
\end {center}

    Note that A is $\gamma$-open [9] iff A =$int_\gamma$(A). A set A is
called $\gamma$- closed [7] iff X-A is $\gamma$-open.\\
{\bf Definition [8].}  A point x $\in$ X is called a
$\gamma$-closure point of A $\subseteq$ X, if $U^\gamma\cap
A=\phi$, for each open nbd. U of x. The set of all
$\gamma$-closure points of A is called $\gamma$-closure of A and
is denoted by $cl_\gamma$(A). A subset A of X is called
$\gamma$-closed, if $cl_\gamma(A)\subseteq$A. Note that
$cl_\gamma$(A) is contained in every $\gamma$-closed superset
of A.\\
{\bf Definition [6].} An operation $\gamma : \tau \rightarrow
P(X)$ is said be $\gamma$-open, if $V^\gamma$ is $\gamma$ open for
each $V\in\tau$.\\
{\bf Definition [9].} An operation $\gamma$ on $\tau$ is said be
regular, if for any open nbds U,V of x$\in$ X, there exists an
open nbd W of x such that $U^\gamma \cap V^\gamma\supseteq
W^\gamma$.\\
{\bf Definition [9].}A mapping $f:(X,\tau_1)\rightarrow(Y,\tau_2)$
is said to be $(\gamma,\beta)$-continuous iff for each $x\in X$
and each open set V containing f(x), there exists an open set U
such that $x\in U$and $f(U^\gamma)\subseteq V^\beta$, where
$\gamma : \tau_1 \rightarrow P(X)$; $\beta : \tau_2\rightarrow
P(Y)$ are operations
on $\tau_1$ and $\tau_2$ respectively.\\
{\bf Definition [2].}A mapping $f:(X,\tau_1)\rightarrow(Y,\tau_2)$
is said to be $(\gamma,\beta)$-closed (respt.
$(\gamma,\beta)$-open), if for any $\gamma$-closed(respt.
$\gamma$-open)set A of X, f(A) is $\beta$-closed
(respt. $\beta$-open) in Y.\\
{\bf Definition [9].} A mapping
$f:(X,\tau_1)\rightarrow(Y,\tau_2)$ is said to be
$(\gamma,\beta)$-homeomorphic, if f is bijective,
$(\gamma,\beta)$-continuous and $f^{-1}$ is $(\beta,\gamma)$-continuous.\\
{\bf Definition [9].} A space X is called a $\gamma$-$T_2$ space,
if for each distinct points x,y $\in$X, there exist open sets U,V
such that $x\in U$ and $y\in V$ and $U^\gamma\cap V^\gamma =\phi$.
\section {$\gamma_0$-compact and $\gamma^*$-regular spaces.}
{\bf Definition 1 [5].} A subset A of a space X is
$\gamma_0$-compact, if every cover $\{V_i:i\in I\}$ of X by
$\gamma$-open sets of X, there exists a finite subset $I_0$ of I
such that $A\subseteq \bigcup_{i\in I_0}cl_\gamma(V_i)$. A space X
is $\gamma_0$-compact [5], if X = $\bigcup_{i \in
I_0}cl_\gamma(V_i)$,
for some finite subset $I_0$ of I.\\
{\bf Definition [9].}  An operation $\gamma$ on $\tau$ is said to
be open if for every nbd U of each $x\in X$,there exists
$\gamma$-open set B such that $x\in B$ and  $U^\gamma\subseteq B$.

    In [5], in the sequel, we shall make use of the following
results:
\\
{\bf Theorem [5].}  Let X be a $\gamma_0$-compact space. Then each
$\gamma$-closed subset of X is $\gamma_0$-compact, where $\gamma$
is a regular operation.
\\
{\bf Theorem [5].} Let $f:(X,\tau_1)\rightarrow(Y,\tau_2)$ be a
bijective $(\gamma,\beta)$-continuous function and $\beta$ be
open. Then f is $(\gamma,\beta)$-open (respt. closed), if and only
if $f^{-1}$ is $(\beta,\gamma)$-continuous.
\\
{\bf Theorem [5].} Let $f:(X,\tau_1)\rightarrow(Y,\tau_2)$ be a
$(\gamma,\beta)$-continuous from a $\gamma_0$-compact space X onto
a space Y. Then Y is $\gamma_0$-compact, where $\beta$ is open.

      We know [5] that G is $\gamma$-open in A X iff G = $A\cap O$, where O
is $\gamma$-open set in X. Using this definition, we prove the
following:
\\
{\bf Theorem 1.} Let $f:(X,\tau_1)\rightarrow(Y,\tau_2)$ be a
$(\gamma,\beta)$-continuous injective function, where $\beta$ is
open . Let C be a $\gamma_0$-compact subspace of a space X. Then
f(C) is $\gamma_0$-compact in Y.
\\
{\bf Proof.} Let $\{V_i:i\in I\}$ be any cover by $\beta$-open
sets of f(C). Then
\begin{eqnarray}
    V_i = U_i \cap f(C)
\end{eqnarray}
where each $U_i$ is $\beta$-open in Y. Since f is
$(\gamma,\beta)$-continuous, therefore $f^{-1}(U_i)$ is
$\gamma$-open in X. Since f is one-one, therefore, (2.1) gives
$f^{-1}(V_i) = f^ {-1}(U_i)\cap f^{-1}f(C)= f^{-1}(U_i)\cap C$.
Clearly $\{W_i:W_i = f^{-1}( V_i)\}$ is $\gamma$-open cover of C.
Since C is $\gamma_0$-compact, therefore there exists a finite
subset $I_0$ of I such that C = $\cup_{i\in
I_0}cl_\gamma(W_i)$.This gives that\\
f(C) = $f(\cup_{i\in I_0}cl_\gamma(W_i))$
 = $\cup_{i\in I_0} fcl_\gamma(W_i)$

    $\subseteq \cup_{i\in I_0}cl_\beta f(W_i)$         ( by Theorem 4.13 [9] )

    = $\cup_{i\in I_0}cl_\beta ff^-1(V_i)$

     = $cl_\beta(V_i)$. Consequently f(C)$\subseteq cl_\beta(V_i)$.\\
This proves that f(C)is $\gamma_0$-compact. This completes the
proof.
\\
{\bf Corollary [5].} Let $f:(X,\tau_1)\rightarrow(Y,\tau_2)$
$(\gamma,\beta)$-continuous from a $\gamma_0$-compact space X onto
a space Y. Then Y is $\gamma_0$-compact, where $\beta$ is open.\\
{\bf Theorem 2.} Let $f:(X,\tau_1)\rightarrow(Y,\tau_2)$ be a
$(\gamma,\beta)$-continuous injective function from a
$\gamma_0$-compact space X into a $\gamma$-$T_2$ space Y and
$\beta$ be open. If $\gamma$ is regular and $\gamma$-open, then f
is
$(\gamma,\beta)$-closed.\\
{\bf Proof.} Let C be a $\gamma$-closed set in X. We show hat f(C)
is $\beta$-closed in Y. Since $\gamma$ is regular, therefore by
theorem 4 [5], C is $\gamma_0$-compact subset of X. Also by
theorem 1, f(C)is $\gamma_0$-compact in Y. Since Y is
$\gamma$-$T_2$ and $\gamma$ is regular and $\gamma$-open,
therefore by theorem 2 [6], f(C) is $\beta$-closed in Y. This
completes the proof.

    The following theorem is immediate:\\
{\bf Theorem 3.} Let $f:(X,\tau_1)\rightarrow(Y,\tau_2)$ be a
$(\gamma,\beta)$-continuous bijective function from a
$\gamma_0$-compact space X onto $\gamma$-$T_2$ space Y. Then f is
$(\gamma,\beta)$-homeomorphism, where $\gamma$ is regular and
$\gamma$-open, and $\beta$ is open.

    The following is easy to prove:\\
{\bf Theorem 4.} If  $A,B\subseteq X$ are $\gamma_0$-compact such
that X = $cl_\gamma(A)\cup cl_\gamma(B)$, then X is
$\gamma_0$-compact, where $\gamma$ is regular and open.\\
{\bf Definition [5].}  A space X is said to be $\gamma^*$-regular
space, if for any $\gamma$-closed set A and $x\in A$, there exist
$\gamma$-open sets U,V such that $x\in U$, $A\subseteq V$ and $U
\cap V = \phi$.

    It is known [5] that $\gamma^*$-regular space is not regular in general. Next
we characterize $\gamma^*$-regular spaces as:\\
{\bf Theorem 5.} A space X is $\gamma^*$-regular iff for each
$x\in X$ and a $\gamma$-closed set A such that $x\notin A$,there
exist $\gamma$-open sets U, V in X such that $x\in U$ and
$A\subseteq V$ and $cl_\gamma(U)\cap cl_\gamma(V) = \phi$.\\
{\bf Proof.} For each $x\in X$ and a $\gamma$-closed set A such
that $x\notin A$, therefore by theorem 9 [4], there is a
$\gamma$-open set W such that
\begin{center}
 $x\in W$, $cl_\gamma(W)\subseteq X- A$.
\end{center}
Again by theorem 9[4], there is a $\gamma$-open set U containing x
such that $cl_\gamma(U)\subseteq W$. Let V = $X-cl_\gamma(W)$.
Then
    $cl_\gamma(U)\subseteq W \subseteq cl_\gamma(W)\subseteq X -
    A$,
implies  $A\subseteq X-cl_\gamma(W) = V$.\\
Also    $cl_\gamma(U)\cap cl_\gamma(V) = cl_\gamma(U)\cap
 cl_\gamma(X-cl_\gamma(W))$
                $\subseteq W \cap cl_\gamma(X-cl_\gamma(W))$
\\
                $\subseteq cl_\gamma(W \cap(X-cl_\gamma(W))$       ( by lemma 2(3)[10])
\\
             = $cl_\gamma(\phi) = \phi$.
\\
Thus U, V are the required $\gamma$-open sets in X. This proves
the necessity. The sufficiency is immediate. This completes the
proof.

\section {$\gamma$-locally compactness.}
{\bf Definition [4].} A $\gamma$-nbd of $x\in X$ is a set U of X
which contains a $\gamma$-open set V containing x. Evidently, U is
a $\gamma$-nbd of x iff $x\in Int_\gamma(U)$.\\
{\bf Definition 2 .} A space X is said to be $\gamma$-locally
compact at $x\in X$ iff x has a $\gamma$-nbd which is
$\gamma_0$-compact in X. If X is $\gamma$-locally compact at every
point, then X is called a $\gamma$-locally compact space.\\
{\bf Example.} Let X= $\{a,b,c\}$, $\tau =\{\phi , X, \{a\},
\{b\},\{a,b\},\{a,c\}\}$. For $b\in X$, define an operation
$\gamma:\tau \rightarrow P(X)$ by
\begin{center}
$\gamma(A) =  \left\{
\begin{array}{ccc}
A, & \mbox{if $b \in A$}\\
cl(A), & \mbox{if $b \not\in A$}\\
\end{array} \right.$\\
\end{center}
Clearly, $\{a,b\},\{a,c\},\{b\}, X ,\phi$ are $\gamma$-open sets
in X. Then

    $U_a = \{\{a,b\},\{a,c\},X\}$

    $U_b = \{\{b\},\{a,b\},X\}$

    $U_c = \{\{a,c\},X\}$
\\
are the $\gamma$-nbd systems at a,b,c. Routine calculations shows
that X is $\gamma$-locally compact at each, $x\in X$. Hence X is
$\gamma$-locally compact.

    $\gamma_0$-compact spaces are $\gamma$-locally compact, since X is
$\gamma$-nbd of its points which is $\gamma_0$-compact. The
converse
is not true in general. Thus we have the following:\\
{\bf Theorem 6.} Every $\gamma_0$-compact space is
$\gamma$-locally compact.

    The following theorem shows that $\gamma$-locally compactness is
$\gamma$-closed hereditary property:\\
{\bf Theorem 7.} Every $\gamma$-closed subspace of
$\gamma$-locally compact space is $\gamma$-locally compact, where
$\gamma$ is regular.\\
{\bf Proof.} Let X be a $\gamma$-locally compact space and A a
$\gamma$-closed set in X. Let $x\in A$. Since X is
$\gamma$-locally compact, therefore for $x\in X$, there is a
$\gamma$-nbd V of $x\in X$ such that V is $\gamma_0$-compact. Let
$U = A \cap V$. Then U is a $\gamma$-nbd of x in A, since V is
$\gamma$-nbd of $x\in X$. Clearly U is $\gamma$-closed in V, since
A is $\gamma$-closed in X. Thus by theorem 4 [5], U is
$\gamma_0$-compact in V. This proves that A is $\gamma$-locally
compact at any point $x\in A$. This completes the proof.\\
{\bf Theorem 8.} If X is a $\gamma$-locally compact and
$\gamma$-$T_2$ space, then for all $x\in X$ and for all $\gamma$-
nbds U of x, there exists a $\gamma$-nbd V of x which is
$\gamma_0$-compact and such that $V\subseteq U$, where $\gamma$ is
regular and open.\\
{\bf Proof.} Let $x\in X$. Since X is $\gamma$-locally compact,
therefore x has a $\gamma$-nbd K which is $\gamma_0$-compact. Let
U be any $\gamma$-nbd of x. Then W = $Int_\gamma(K \cap U)$ is
$\gamma$-open in X. $W \subseteq K$ implies $cl_\gamma(W)
\subseteq cl_\gamma(K)$. Then by theorem 2 [6], $cl_\gamma(K)$ = K
gives $cl_\gamma(W)\subseteq K$. Also by theorem 4 [5],
$cl_\gamma(W)$ is $\gamma_0$-compact in X. Put $C = cl_\gamma(W) -
W$. Clearly, C is $\gamma$-closed and hence $\gamma_0$-compact in
$cl_\gamma (W)$. Moreover $x\notin C$, then by theorem 1 [6],
there exist open sets G, H in $cl_\gamma(W)$ such that $x\in
G^\gamma$, $C\subseteq H^\gamma$ and $G^\gamma \cap H^\gamma =
\phi$. Thus G is a $\gamma$-nbd of x in $cl_\gamma(W)$ and
$cl_\gamma(W)$ is a $\gamma$-nbd of x in X imply G is a
$\gamma$-nbd of x in X.
 $G^\gamma \cap H^\gamma = \phi$ implies
\begin{eqnarray}
G^\gamma \subseteq cl_\gamma(W) - H^\gamma \subseteq cl_\gamma(W)
- C = W
\end{eqnarray}
$cl_\gamma(W)- H^\gamma$ is $\gamma$-closed in $cl_\gamma(W)$ and
so is in X. By (3.1),
\begin{center}
$G^\gamma \subseteq cl_\gamma(G^\gamma)\subseteq
cl_\gamma(cl_\gamma(W)- H^\gamma) = cl_\gamma(W)- H^\gamma
\subseteq W \subseteq U$.
\end{center}
 Put $cl_\gamma(G^\gamma) = V$. Then by theorem 4 [5], V
is $\gamma_0$-compact in $cl_\gamma(W)$ and hence in X.
 Consequently, $x\in V \subseteq U$, where V is $\gamma_0$-compact
$\gamma$-nbd of x. This completes the proof.

    It is known [6], that in $\gamma$-$T_2$ space, every $\gamma_0$-compact
subset of $\gamma$-$T_2$ space is $\gamma$-closed, where $\gamma$
is regular and $\gamma$-open. We use this fact and prove the
following:\\
{\bf Theorem 9.} If a space X is $\gamma$-locally compact, then
each of its points is a $\gamma$-interior point of some
$\gamma_0$-compact subspace of X.\\
{\bf Proof.} Since X is $\gamma$-locally compact, then each $x\in
X$ has a $\gamma$-nbd N such that $cl_\gamma$(N)is
$\gamma_0$-compact, that is, $cl_\gamma(N)$ is a
$\gamma_0$-compact $\gamma$-nbd of x and so x is a
$\gamma$-interior point. Thus each $x\in X$ is a $\gamma$-interior
point of some $\gamma_0$-compact subspace of X. This completes the
proof.

    The partial converse of theorem 9 is:\\
{\bf Theorem 10.} If in a $\gamma$-$T_2$ space X, each $x\in X$ is
a $\gamma$-interior point of some $\gamma_0$-compact subspace of
X, then X is $\gamma$-locally compact, where $\gamma$ is regular
and
$\gamma$-open.\\
{\bf Proof.} Let $x\in X$. Then by hypothesis, there exists a
$\gamma_0$-compact subspace C of X such that x is the
$\gamma$-interior point of C. Since X is $\gamma$-$T_2$, by
theorem 2 [6], C is $\gamma$-closed subset of X. $cl_\gamma(C)$ =
C gives that every $x\in X$ has a $\gamma$-nbd C whose
$\gamma$-closure is $\gamma_0$-compact. This proves that X is
$\gamma$-locally compact.
 This completes the proof.

    Combining theorems 9 and 10, we have:\\
{\bf Theorem 11.} A $\gamma$-$T_2$ space X is $\gamma$-locally
compact iff each $x\in X$ is a $\gamma$-interior point of some
$\gamma_0$-compact subspace of X, where $\gamma$ is regular and
$\gamma$-open.\\
{\bf Definition 3.} A property is called $(\gamma,\beta)$-
topological property, if it is possessed by a space X, then it is
also possessed by all spaces $(\gamma,\beta)$-homeomorphic to X.

    In [5], it is proved that if a surjective function
$f:(X,\tau_1)\rightarrow(Y,\tau_2)$ is $(\gamma,\beta)$-continuous
and X is a $\gamma_0$-compact, then Y is $\gamma_0$-compact, where
$\beta$ is open.

    We use this result and prove the following theorem which shows that
$\gamma$-locally compactness is $(\gamma,\beta)$-topological
property.\\
{\bf Theorem 12.} Let $f:(X,\tau_1)\rightarrow(Y,\tau_2)$ be
$(\gamma,\beta)$-open,$(\gamma,\beta)$-continuous surjection. If X
is $\gamma$-locally compact, then Y is also $\gamma$-locally
compact, where $\beta$ is open.\\
{\bf Proof.} Let $y\in Y$. Then there is $x\in X$ such that f(x) =
y. Let $U_x$ be a $\gamma_0$-compact,$\gamma$-nbd of x, since X is
$\gamma$-locally compact. Also $f(x)\in
f(Int_\gamma(U_x))\subseteq f(U_x)$ gives $f(U_x)$ is $\gamma$-nbd
of f(x). Since f is $(\gamma,\beta)$-open, then
$f(Int_\gamma(U_x))$ is $\gamma$-open in Y. Also f is
$(\gamma,\beta)$-continuous, therefore by theorem 3 [5], $f(U_x)$
is $\gamma_0$-compact. Thus each $y\in Y$ has a
$\gamma_0$-compact, $\gamma$-nbd in Y. This proves that Y is
$\gamma$-locally compact. Hence the proof.

    In [5], it is proved that a $\gamma$-closed subset of a
$\gamma_0$-compact spaces is $\gamma_0$-compact, if $\gamma$ is
regular. We use this result and prove that $\gamma$-locally
compactness is a $\gamma$-closed hereditary property.\\
{\bf Theorem 13.} Every $\gamma$-closed subspace of a
$\gamma$-locally compact space is $\gamma$-locally compact, if
$\gamma$ is regular.\\
{\bf Proof.} Let X be a $\gamma$-locally compact space and A a
$\gamma$-closed set in X. Let $a\in A$. Since X is
$\gamma$-locally compact, therefore for $a\in X$, there is a
$\gamma$-nbd V of a such that V is $\gamma_0$-compact in X. Let $U
= V \cap A$. Then U is a $\gamma$-nbd of A. Clearly U is
$\gamma$-closed in V and so by theorem 4 [5], U is
$\gamma_0$-compact in A, that is, A is $\gamma$-locally compact at
a. This proves that A is $\gamma$-locally compact. Hence the
theorem.
\section {$\gamma$-normal spaces.}

    In [5], $\gamma$-normal spaces have been defined and
studies. Here we prove another interesting characterization of
such spaces.\\
{\bf Definition [5].} A space X is said to be $\gamma$-normal
space, if for any disjoint $\gamma$-closed sets A, B of X, there
exist $\gamma$-open sets U, V such that $A \subseteq U$, $B
\subseteq V$ and $U \cap V = \phi$.

    The following examples show that $\gamma$-normality and normality
are independent notions:\\
{\bf Example 1.} Let X= $\{a,b,c\}$, $\tau =\{\phi , X, \{a\},
\{b\},\{a,b\},\{a,c\}\}$. For $b\in X$, define an operation
$\gamma:\tau \rightarrow P(X)$ by
\begin{center}
$\gamma(A) =  \left\{
\begin{array}{ccc}
A, & \mbox{if $b \in A$}\\
cl(A), & \mbox{if $b \not\in A$}\\
\end{array} \right.$\\
\end{center}
Then X is $\gamma$-normal. X is not normal, since for closed sets
$\{a,c\},\{d\}$, there do not exist disjoint open sets containing
$\{a,c\}$
and $\{d\}$ respectively.\\
{\bf Example 2.} Let X= $\{a,b,c\}$, $\tau =\{\phi , X, \{a\},
\{b\},\{a,b\},\{a,c\}\}$. For $b\in X$, define an operation
$\gamma:\tau \rightarrow P(X)$ by
\begin{center}
$\gamma(A) =  \left\{
\begin{array}{ccc}
A, & \mbox{if $b \in A$}\\
cl(A), & \mbox{if $b \not\in A$}\\
\end{array} \right.$\\
\end{center}
Then X is not $\gamma$-normal, since for $\gamma$-closed sets
$\{a\},\{c\}$, there do not exist disjoint $\gamma$-open sets
containing $\{a\},\{c\}$ respectively, where X is normal.

    It is known [10],that if A is  $\gamma$-open,then $A \cap cl_\gamma(B)
\subseteq cl_\gamma(A \cap B)$  ...      (3.2).\\
In [5], $\gamma$-normal spaces have been characterized as :\\
{\bf Theorem A [5].} A space X is $\gamma$-normal iff for any
$\gamma$-closed set A and $\gamma$-open set U containing A, there
is a $\gamma$-open set V containing A such that $cl_\gamma(V)
\subseteq U.$

    We use (3.2) and theorem A, and prove the following:\\
{\bf Theorem 14.} A space X is $\gamma$-normal iff for each pair
A, B of disjoint $\gamma$-closed sets in X, there exist
$\gamma$-open sets U, V in X such that $A \subseteq U$, $B
\subseteq V$ and $cl_\gamma(U)\cap cl_\gamma(V)= \phi$.\\
{\bf Proof.} The sufficiency is clear. We prove only the
necessity. Let A be a $\gamma$-closed set and B be a
$\gamma$-closed set not containing A. Then X - B is $\gamma$-open
and $A \subseteq X - B$. Then by theorem A [5], there is a
$\gamma$-open set C such that
\begin{center}
     $A \subseteq C \subseteq cl_\gamma(C)\subseteq X- B.$
\end{center}
Since $cl_\gamma(C)\subseteq X- B$, again by theorem A [5], there
is a $\gamma$-open set U containing A such that
$cl_\gamma(U)\subseteq C $. Consequently,
\begin{center}
 $A \subseteq U \subseteq cl_\gamma(U)\subseteq C \subseteq cl_\gamma(C) \subseteq X -
 B.$
\end{center}
$cl_\gamma(C)\subseteq X - B$ implies $B \subseteq X - cl_\gamma
(C)$. Put V = X - $cl_\gamma(C)$. Then V is $\gamma$-open
containing
B and moreover\\
    $cl_\gamma(U)\cap cl_\gamma(V) = cl_\gamma(U)\cap cl_\gamma(X -
    cl_\gamma(C))$
                     $\subseteq C \cap cl_\gamma(X -
                     cl_\gamma(C))$ \\
                      $\subseteq cl_\gamma(C \cap(X - cl_\gamma(C))$         ( by
                      (3.2))\\
                   $= cl_\gamma(\phi) = \phi$ .\\
Thus U, V are the required $\gamma$-open sets in X. This proves
the necessity. Hence the theorem.

\end{document}